\documentclass[leqno, 10pt]{amsart}
\usepackage{amsmath, amsthm, amssymb, latexsym}
 \newtheorem{definition}{Definition}[section]

 \newtheorem{proposition}[definition]{Proposition}

 \newtheorem*{theorem*}{Theorem}
\newtheorem*{proposition*}{Proposition}
\newtheorem*{lemma*}{Lemma}

 \theoremstyle{remark}

  \newtheorem*{acknowledgements}{Acknowledgements}

\newcommand{\op}[1]{\operatorname{#1}}

\newcommand{\acou}[2]{\ensuremath{\langle #1 , #2 \rangle}}

\newcommand{\tr}{\op{tr}}
\newcommand{\Tra}{\ensuremath{\op{Trace}}}

\newcommand{\TR}{\ensuremath{\op{TR}}}

\newcommand{\Res}{\ensuremath{\op{Res}}}
\newcommand{\res}{\ensuremath{\op{Res}}}

\newcommand{\bint}{\ensuremath{-\hspace{-2,4ex}\int}}

\newcommand{\C}{\ensuremath{\mathbb{C}}}

\newcommand{\Z}{\ensuremath{\mathbb{Z}}} 
\newcommand{\CZ}{\ensuremath{\mathbb{C}\!\setminus\!\mathbb{Z}}}

\newcommand{\Ca}[1]{\ensuremath{\mathcal{#1}}}

\newcommand{\cE}{\Ca{E}}

\newcommand{\cH}{\ensuremath{\mathcal{H}}}

\newcommand{\cR}{\ensuremath{\mathcal{R}}}
\newcommand{\cS}{\ensuremath{\mathcal{S}}}

\newcommand{\sD}{\ensuremath{{/\!\!\!\!D}}}%{/{\!\!\!\!D}}
\newcommand{\sS}{\ensuremath{{/\!\!\!\!\!\;S}}}

\newcommand{\pdoz}{\ensuremath{\Psi^{\Z}}} 
\newcommand{\pdocz}{\ensuremath{\Psi^{\CZ}}} 
 
\newcommand{\psido}{$\Psi$DO} 
\newcommand{\psidos}{$\Psi$DOs}

\newcommand{\ord}{{\op{ord}}}

\newcommand{\End}{\ensuremath{\op{End}}}

\newcommand{\rk}{\op{rk}}

\newcommand{\Vol}{\op{Vol}}

\begin{document}
\title{Noncommutative geometry and lower dimensional volumes in Riemannian geometry} 

\author{Rapha\"el Ponge}

\address{Department of Mathematics, University of Toronto, Canada.}
% \address{Department of Mathematics, University of Toronto, 40 St George Street, Toronto, ON M5S 2E4, Canada.}
\email{ponge@math.toronto.edu}
\thanks{Research partially supported by NSERC grant 341328-07  and 
 by a new staff matching grant from the Connaught Fund of the University of Toronto}
\begin{abstract}
 In this paper we explain how to define ``lower dimensional'' volumes of any compact Riemannian manifold as the integrals of local Riemannian 
 invariants. For instance we give sense to the area and the length of such a manifold in \emph{any} dimension.  
 Our reasoning is motivated by an idea of Connes and involves in an essential way noncommutative geometry and the analysis of Dirac operators on spin 
 manifolds. However, the ultimate definitions of the lower dimensional volumes don't involve noncommutative geometry or spin structures at all.
\end{abstract}

\maketitle 
\numberwithin{equation}{section}

 \maketitle 
    
\section*{Introduction}
An important application of Riemannian geometry is to allow us to define the volume element of a Riemannian manifold $(M^{n},g)$ as the  
volume form (or $1$-density) $dv_{g}(x)=\sqrt{g(x)}d^{n}x$, in such way that the volume of an open set $U\subset M$ is given by integrating $dv_{g}(x)$ over $U$. 

It would be interesting to push further this idea and to ask for the existence of a \emph{length element} $ds$ in such way to be able to define the 
$k$'th power $ds^{k}$ for any $k=2,\ldots,n$ so that for $k=n$ we recover the volume element $dv_{g}(x)$. Then $ds^{k}$ could be interpreted as 
the $k$'th dimensional volume element, e.g.~$ds^{2}$ would be the area element of $M$. Furthermore, if for any integer $k\leq n-1$ we were able find 
a way to integrate $ds^{k}$,  then we would be able define the $k$'th dimensional volume of $(M,g)$ as the integral of $ds^{k}$. 

At first glance it seems difficult to work out out the above ideas by remaining within the framework of classical Riemannian geometry. 
The equation $ds^{n}=dv_{g}(x)$ suggests that 
$ds$ should be a $n$'th root of the volume element $dv_{g}(x)$ in a sense which has yet to be defined. Furthermore, if we see $dv_{g}(x)$ as a 
$1$-density, then it would be natural to think that $ds$ could be a $\frac{1}{n}$-density. However, the classical integral only allows us to 
integrate $1$-densities, so we wouldn't be able to use the classical integral to integrate $ds^{k}$ for $k \leq n-1$. 

The aim of this note is to propose to make use of noncommutative geometry to resolve the above issues. Our approach is motivated by 
an idea of Connes~(\cite{Co:GCMFNCG}, \cite{CC:SAP}) in connection with his spectral 
interpretation of the Einstein-Hilbert action. The latter is given by the integral of the scalar curvature and yields 
the contribution of gravity forces to the functional action of the Standard Model. Connes observed that when $M$ is spin the Einstein-Hilbert action is a 
constant multiple of the noncommutative residue $\bint \sD^{-n+2}$, where $\sD$ denotes the Dirac operator on $M$ (see~\cite{KW:GNGWR}, \cite{Ka:DOG}). 
This result is an important ingredient in 
the spectral action principle of~\cite{CC:SAP} and its recent refinement in~\cite{CC:GSMNM}. This also plays an important role in the axiomatic 
characterization of spin manifolds (see~\cite{Co:GCMFNCG}, \cite{RV:RMBCG}). 

The noncommutative residue of Wodzicki~(\cite{Wo:LISA}, \cite{Wo:NCRF}) and Guillemin~\cite{Gu:NPWF} is a trace on the algebra of 
(integer order) \psidos\ on $M$. An important feature is that 
it allows us to extend to all \psidos\ the Dixmier trace, which plays the role of the integral in the framework of noncommutative geometry. Thus it 
allows us to define the (noncommutative) integral of a \psido\ even if the latter is not in the domain of the Dixmier trace. In QFT the inverse of the inverse of the 
Dirac operator 
is the free propagator for Fermions, so it has the dimension of a length. This motivated Connes to interpret in dimension $4$ the noncommutative residue 
$\bint \sD^{-4+2}=\bint \sD^{-2}$ as the area in Planck's units of the manifold $M$. 

In this paper we push further Connes' idea to define lower dimensional dimensional volumes in Riemannian geometry. The outcome 
depends on the parity of $n$. Let us first assume $n$ even and $M$ spin. Our starting point is the well known fact that the 
noncommutative residue density of the operator $\sD^{-n}$ allows us to recapture the volume form of $M$. We thus can interpret $\sD^{-n}$ as the 
 \emph{noncommutative volume element} of $M$. The upshot of this is twofold:\smallskip  

(i) While there are some issues with defining the $n$'th root of a volume form, the $n$'th root a positive operator makes well sense. We thus 
can define the  \emph{noncommutative length element} $ds$ as a constant multiple of $|\sD|^{-1}$.\smallskip  

(ii) Whereas the classical integral allows us to integrate $1$-densities only,  the noncommutative residue enables us to 
integrate any \psido. Therefore, for $k=1,\ldots,n$ we can define the \emph{$k$'th dimensional volume}  $ \Vol^{(k)}_{g}M$ of $(M,g)$ by letting 
\[
 \Vol^{(k)}_{g}M=\bint ds^{k}.
\]
\indent The next step is to express the lower dimensional volumes $ \Vol^{(k)}_{g}M$ in a purely differential geometric way. 
To this end we use the well known connection between the noncommutative 
residues $\bint |\sD|^{-k}$ and the coefficients of the heat kernel asymptotics of $\sD^{2}$. These coefficients are fairly well 
understood thanks to the work of Atiyah-Bott-Patodi, Gilkey and others. In particular, they are shown to be local Riemannian 
invariants, i.e., they can be written a universal linear combination in complete tensorial contractions of covariant derivatives of the curvature 
tensor. This allows us to shows that $\Vol^{(k)}_{g}M$ is the integral of a local Riemannian invariant, which vanishes when $k$ is odd 
(see Proposition~\ref{prop:Even.main} for the precise statement).

Now, the local Riemannian invariants makes sense independent of the existence of a spin structure, we can use the geometric expression for 
$\Vol^{(k)}_{g}M$ to extend its definition to general, possibly non-spin, compact Riemannian dimension (Definition~\ref{def:Even.general-def}). 
We thus get a purely Riemannian definition of the lower dimensional volumes without 
reference anymore to noncommutative geometry. In particular, we give explicit formulas for the area $\op{Area}_{g}M:=\Vol_{g}^{(2)}M$ 
in dimension 4 and in dimension 6. In dimension 4 the formula agrees with Connes' formula for the area of a 4-dimensional spin manifold. 

In odd dimension the approach is similar. We first $M$ spin and we define the noncommutative length element $ds$ as a constant multiple of 
$|\sD|^{-1}$. For $k=1,\ldots,n$ we define $\Vol_{g}^{(k)}M$ to be $\bint ds^{k}$ and we express it  as the integral of a local Riemannian invariant 
(see Proposition~\ref{prop:Odd.main}). In contrast with the even dimensional case $\Vol_{g}^{(k)}M$ vanishes when $k$ even. We then use 
these formulas to define the lower dimensional volumes on any compact Riemannian manifold of odd dimension (Definition~\ref{def:Odd.general-def}). 
As examples we work out explicit formulas for the length $\op{Length}_{g}M:=\Vol^{(1)}_{g}M$ in dimension 3 and in dimension~5. 

This paper is organized as follows. In Section~\ref{sec:NCG} we recall some important fact about Connes's quantized calculus and the noncommutative 
residue.  In Section~\ref{sec:even-dimension} we define the lower dimensional volumes in even dimension. Finally, in Section~\ref{sec:odd-dimension}
 we deal with the lower  dimensional volumes in odd dimension.

\begin{acknowledgements}
 I wish to thank for its hospitality the University of California at Berkeley where this paper reached its final form.
\end{acknowledgements}

\section{Quantized calculus and the noncommutative residue}
\label{sec:NCG}
In this section we recall some definitions and properties of Connes' quantized calculus and the noncommutative residue trace of Wodzicki and Guillemin. 

The quantized calculus of Connes~\cite{Co:NCG} provides us with a dictionary for translating tools of 
the infinitesimal calculus into the language of quantum mechanics. Given a separable Hilbert space $\cH$ the first few lines of this dictionary are: 

\begin{center}
    \begin{tabular}{c|c}  
        Classical & Quantum \\ \hline
        
      Real variable &  Selfadjoint operator on $\cH $ \\  
    %    & & \\ 
       Complex variable & Operator on $\cH $ \\
 Infinitesimal variable & Compact operator on $\cH $ \\
    %    & & \\ 
       Infinitesimal of order  & Compact operator $T$ such that \\ 
                $\alpha>0$       &  $\mu_{k}(T)=\op{O}(k^{-\alpha})$\\
                 Integral & Dixmier Trace $\bint$ 
    \end{tabular}
\end{center}

Here $\mu_{n}(T)$ denotes the  $(n+1)$'th characteristic value of $T$, i.e., the $(n+1)$'th eigenvalue 
of $|T|=(T^{*}T)^{\frac12}$. In particular, by the min-max principle  we have 
\begin{equation}
    \mu_{n}(T)  =  \inf\{ \|T_{E^\perp}\|; \dim E=n\}.
%     \nonumber \\
\end{equation}

The Dixmier trace is defined on infinitesimal operators of order~$\leq 1$ and arises in the analysis of the logarithmic divergency of the partial 
sums $\sigma_{N}(T):=\sum_{k<N}\mu_{k}(T)$. For instance, for a compact operator $T\geq 0$ we have
\begin{equation}
    \lim_{N \rightarrow \infty}\frac{1}{\log N}\sum_{k<N}\mu_{k}(T) \Longrightarrow \bint T=L.
%     \label{¥}
\end{equation}
In general, if $\mu_{k}(T)=\op{O}(k^{-1})$ then $\sigma_{N}(T)=\op{O}(\log N)$ and we can extract a limit point of the sequence 
$(\frac{\sigma_{N}(T)}{\log N})_{N\geq 2}$ in such way to get a trace (see~\cite{Di:ETNN}, \cite[Appendix A]{CM:LIFNCG}). We then say that $T$ is 
\emph{measurable} when the value of the limit point is 
independent of the limit process and we then denote $\bint T$ this value: this is the \emph{Dixmier trace} of $T$. 

Now, let $M^{n}$ be a compact Riemannian manifold and let $\cE$ be a Hermitian vector bundle over $M$. For $m \in \C$ we let $\Psi^{m}(M,\cE)$ be the 
class of \psidos\ of order $m$ on $M$ acting on the sections of $\cE$. Thus an element of $\Psi^{m}(M,\cE)$ is a continuous operator 
$P:C^{\infty}(M,\cE)\rightarrow C^{\infty}(M,\cE)$ whose Schwartz kernel is smooth off the diagonal and in local coordinates $P$ is of the form
\begin{equation}
    Pu(x)=(2\pi)^{-n}\int_{M}e^{ix.\xi}p(x,\xi)\hat{u}(\xi)d\xi +Ru(x),
%     \label
\end{equation}
where $p(x,\xi)\sim \sum_{j\geq 0} p_{m-j}(x,\xi)$ is a polyhomogeneous symbol of degree $m$ and $R$ is a smoothing operator. 

If $\Re m\leq 0$ then $P$ is bounded on the Hilbert space $L^{2}(M,\cE)$ and is actually a compact operator if we further have $\Re m<0$. In fact, we have: 

\begin{proposition}[Connes~\cite{Co:AFNCG}]\label{prop:QC.Dixmier-trace} 
    Let $P\in \Psi^{m}(M,\cE)$, $\Re m<0$. Then:\smallskip
    
    1) $P$ is an infinitesimal operator of order~$\leq \frac{|\Re m|}{n}$.\smallskip
    
    2) If $\ord P=-n$ then $P$ is measurable for the Dixmier trace and we have 
    \begin{equation}
        \bint P= \frac{1}{n} \Res P,
        \label{eq:QC.Dixmier-NCR}
    \end{equation}
    where $\Res $ denotes the noncommutative residue trace of Wodzicki~(\cite{Wo:LISA}, \cite{Wo:NCRF}) and Guillemin~\cite{Gu:NPWF}. 
\end{proposition}

Before recalling the main definitions and properties of the noncommutative residue trace, let us first notice that Eq.~(\ref{eq:QC.Dixmier-NCR}) 
allows to  extend the Dixmier trace to the whole algebra $\pdoz(M,\cE)$  of integer order \psidos\ by letting
\begin{equation}
    \bint P= \frac{1}{n}\Res P \qquad \text{for any $P\in \pdoz(M,\cE)$}.
    \label{eq:QC.extension-Dixmier-trace}
\end{equation}
In other words the noncommutative residue allows us to integrate any operator $P\in \pdoz(M,\cE)$ even when $P$ is not an infinitesimal of order~$\leq 1$.

The noncommutative residue is a trace on the algebra $\pdoz(M,\cE)$ independently found by Wodzicki~(\cite{Wo:LISA}, \cite{Wo:NCRF})
and Guillemin~\cite{Gu:NPWF}. It can be defined as follows. 

Let $P\in \Psi^{m}(M,\cE)$, $m\in \Z$. Then in local coordinates the Schwartz kernel $k_{P}(x,y)$ has a behavior near the diagonal $y=x$ of the 
form
\begin{equation}
    k_{P}(x,y)=\sum_{-(m+n)\leq j\leq -1}a_{j}(x,x-y)-c_{P}(x)\log|x-y|+\op{O}(1),
    \label{eq:NCR.log-sing}
\end{equation}
with $a_{j}(x,y)$ is homogeneous in $y$ of degree $j$ and $c_{P}(x)$ given by: 
\begin{equation}
    c_{P}(x)=(2\pi)^{-n}\int_{S^{n-1}}p_{-n}(x,\xi)d^{n-1}\xi,
%     \label{¥}
\end{equation}
where $p_{-n}(x,\xi)$ is the symbol of degree $-n$ of $P$. 

As observed by Connes-Moscovici~\cite{CM:LIFNCG} the coefficient $c_{P}(x)$ makes sense globally on $M$ as an $\End\cE$-valued 1-density (i.e.~as a 
section of $|\Lambda|(M)\otimes \End \cE$ where $|\Lambda|(M)$ is the bundle of $1$-densities over $M$).
We then define the \emph{noncommutative residue} of $P$ by means of the formula
\begin{equation}
    \Res P=\int_{M}\tr_{\cE}c_{P}(x).
     \label{eq:NCR.formula-cP}
\end{equation}

Next, the density $c_{P}(x)$ is intimately related to the local analytic extension of the usual trace to \psidos\ of non-integer orders. If $P$ is a 
\psido\ of order $m$ with $\Re m<-n$, then the restriction to the diagonal $k_{P}(x,x)$ of its Schwartz kernel defines a smooth density.  Therefore $P$ is 
trace-class and we have 
\begin{equation}
    \Tra P= \int_{M}\tr_{\cE} k_{P}(x,x).
%     \label
\end{equation} 

In fact the map $P\rightarrow k_{P}(x,x)$  can be analytically extend into a map $P\rightarrow t_{P}(x)$ defined on the class $\pdocz(M,\cE)$ of 
non-integer order \psidos. Here analyticity is meant with respect to holomorphic families of \psidos\ as in~\cite{Gu:RTCAFIO} and~\cite{KV:GDEO}. 
We thus obtain the unique analytic extension of the usual trace to $\pdocz(M,\cE)$ by letting
\begin{equation}
    \TR P =\int_{M}\tr_{\cE}t_{P}(x) \qquad \text{for any $P\in \pdocz(M,\cE)$}.  
%     \label
\end{equation}

Furthermore, if $P$ is a \psido\ of 
integer order $m$ and if $(P(z))_{z\in \C}$ is a holomorphic family of \psidos\ such that $P(0)=P$ and $\ord P(z)=z+\ord P$, then the map 
$z\rightarrow t_{P(z)}(x)$ has at worst a simple pole singularity near $z=0$ such that: 
\begin{equation}
    \Res_{z=0}t_{P(z)}(x)=-c_{P}(x).
     \label{eq:NCR.residue}
\end{equation}
In particular we see that $\Res P=-\Res_{z=0} \TR P(z)$. 

All this shows that the noncommutative residue $\Res$ is (up to a sign) 
the residual functional induced on $\pdoz(M,\cE)$ by the analytic extension to $\pdocz(M,\cE)$ of the 
usual trace. Granted this it is not difficult  to see that the noncommutative residue is a trace on the algebra $\pdoz(M,\cE)$. In fact, as proved by 
Wodzicki~\cite{Wo:PhD} (see also~\cite{Gu:RTCAFIO}) 
the noncommutative residue even is the unique trace on $\pdocz(M,\cE)$ up to constant multiple 
when $M$ is connected. 

Finally, let $P$ be a positive elliptic differential operator of order $m$ with a positive principal symbol and for $t>0$ let $k_{t}(x,y)$ be 
the Schwartz kernel of $e^{-tP}$. Then as $t \rightarrow 0^{+}$ we have an asymptotics in $C^{\infty}(M,|\Lambda|(M)\otimes \End \cE)$ of the form 
\begin{equation}
    k_{t}(x,x) \sim \sum_{j \geq 0} t^{\frac{2j-n}{m}}a_{j}(P)(x),
    \label{eq:NCR.heat-kernel-asymptotics}
\end{equation}
where the coefficients can be locally expressed in terms of the symbol of $P$ (see, e.g., \cite{Gr:AEHE}, \cite{Gi:ITHEASIT}). 

For $\Re s>0$ we have 
$P^{-s}=\Gamma(s)^{-1}\int_{0}^{\infty}t^{s-1}(1-\Pi_{0}(P))e^{-tP}dt$,
where $\Pi_{0}(P)$ denotes the orthogonal projection onto the kernel of $P$. It is standard to use this formula for relating the coefficients of the 
heat kernel asymptotics~(\ref{eq:NCR.heat-kernel-asymptotics}) to the singularities of the zeta function $z\rightarrow t_{P^{-z}}(x)$ (see, e.g., 
\cite{Gi:ITHEASIT}, \cite{Wo:NCRF}). 
In particular, by combining this with~(\ref{eq:NCR.residue}) we see that if for $j=0,1,\ldots,n-1$ we set 
$s_{j}=\frac{n-2j}{m}$, then we have:
\begin{gather}
    mc_{P^{-\frac{n-j}{m}}}(x)=\res_{s=\frac{n-j}{m}}t_{P^{-s}}(x)=0 \qquad \text{if $j$ is odd}, \label{eq:NCR.NCR-heat-kernel-asymptotics1}\\
     mc_{P^{-\frac{n-j}{m}}}(x)=\res_{s=\frac{n-j}{m}}t_{P^{-s}}(x)=\Gamma(\frac{n-j}{m})^{-1}a_{\frac{j}{2}}(P)(x) \qquad \text{if $j$ is even}.
    \label{eq:NCR.NCR-heat-kernel-asymptotics2} 
\end{gather}    
This provides us with an alternative way to compute the densities $c_{P^{-s_{j}}}(x)$. %which we are going to make extensive use in the sequel. 

\section{Lower dimensional volumes (even dimension)}
\label{sec:even-dimension}
In this section we  define the  lower dimensional volumes of a compact Riemannian manifold $(M^{n},g)$ of even dimension. 

Let us first assume $M$ spin and let $\sD$ be its Dirac operator acting on the section of the spin bundle $\sS=\sS^{+}\oplus \sS^{-}$. 
Notice that $|\sD|^{-n}$ allows us to recapture the volume form $dv_{g}(x)=\sqrt{g(x)}d^{n}x$ of $M$. 
Indeed, since $\sD^{2}$ is a Laplace type operator
the  principal symbol of $\sD^{-n}=(\sD^{2})^{-\frac{n}{2}}$ is equal to $(|\xi|_{g}^{2})^{-\frac{n}{2}}=|\xi|_{g}^{-n}$. The latter yields the 
symbol of degree $-n$ of $\sD^{-n}$, so if we work in normal coordinates centered at $x$, then from~(\ref{eq:NCR.formula-cP})  we get
\begin{equation}
    \frac{1}{n}\tr_{\sS} c_{\sD^{-n}}(x)=\frac{(2\pi)^{-n}}{n}\rk \sS.|S^{n-1}|dv_{g}(x)= 
    \frac{(2\pi)^{-\frac{n}{2}}}{\Gamma(\frac{n}{2}+1)}dv_{g}(x) .
    \label{eq:Even.cD-n-volume-form}
\end{equation}
where we have used the fact that $\rk \sS=2^{\frac{n}{2}}$ and 
$|S^{n-1}|=\frac{2\pi^{\frac{n}{2}}}{\Gamma(\frac{n}{2})}=\frac{n\pi^{\frac{n}{2}}}{\Gamma(\frac{n}{2}+1)}$. 
Therefore, for any $f \in C^{\infty}(M)$, we have 
\begin{equation}
    \bint f \sD^{-n}=\frac{1}{n}\int_{M} f(x)\tr_{\sS}c_{\sD^{-n}}(x)= 
   \frac{(2\pi)^{-\frac{n}{2}}}{\Gamma(\frac{n}{2}+1)} \int_{M} f(x)dv_{g}(x) .
     \label{eq:NCG.volume-form-|D|-even}
\end{equation}
Thus the operator $ (2\pi)^{\frac{n}{2}}\Gamma(\frac{n}{2}+1) \sD^{-n}$ recaptures the volume form of $M$. 

Following the idea described in Introduction we seek to define a length element as a $n$'th root of the volume element. Since  
$ (2\pi)^{\frac{n}{2}}\Gamma(\frac{n}{2}+1) \sD^{-n}$ recaptures the volume form, we can interpret it as a (noncommutative) volume element. 
Its $n$'th root then has an obvious meaning, so we are lead to: %the following definition: 

\begin{definition}
 The  noncommutative length element of $M$ is 
\begin{equation}
    ds:=c_{n}|\sD|^{-1}, \qquad c_{n}=\sqrt{2\pi}\Gamma(\frac{n}{2}+1)^{\frac{1}{n}}.
     \label{eq:NCG.length-element-even}
\end{equation}
\end{definition}

Next, the idea to defining for $k=1,\ldots,n$ the $k$'th dimensional is to integrate $ds^{k}$. Notice that $ds^{k}$ is a \psido\ of order $-2k$. 
Since the extension of the Dixmier trace provided by the 
noncommutative residue allows to integrate \emph{any} \psido\ we obtain the definition below.

\begin{definition}
    For $k=1,\ldots,n$ the $k$'th dimensional volume of $M$ is
\begin{equation}
    \Vol^{(k)}_{g}M:= \bint ds^{k}.
     \label{eq:NCG.k-dim-volume.spin}
\end{equation}
\end{definition}

Using~(\ref{eq:NCG.volume-form-|D|-even}) we see that for $k=n$  we have
\begin{equation}
     \Vol^{(n)}_{g}M= (2\pi)^{-\frac{n}{2}}\Gamma(\frac{n}{2}+1) \bint \sD^{-n} = \int_{M}dv_{g}(x)=\Vol_{g}M.
     \label{eq:Even.volume-k=n}
\end{equation}
Hence $\Vol^{(n)}_{g}M$ agrees with the usual volume of $M$. 

We shall now give a differential-geometric interpretation of $\Vol^{(k)}_{g}M$ for $k \leq n-1$  similar to that  for $k=n$.  

In the sequel we let $R_{ijkl}=\acou{R(\partial_{i},\partial_{j})\partial_{k}}{\partial_{l}}$ be the components of the curvature tensor in local 
coordinates. As usual we will use the metric $g=(g_{ij})$ and its inverse  $g^{-1}=(g^{ij})$ to lower and raise indices. For instance the Ricci tensor is 
$\rho_{jk}=R^{i}_{~jki}$ and the scalar curvature is $\kappa=\rho^{j}_{~j}$.  

In terms of the densities $c_{|\sD|^{-k}}(x)$ we have
\begin{equation}
    \Vol^{(k)}_{g}M= (c_{n})^{k}\bint |\sD|^{-k}=  \frac{(c_{n})^{k}}{n}\int_{M} \tr_{\sS}c_{|\sD|^{-k}}(x).
     \label{eq:Even.Volk-cD-k}
\end{equation}
Notice that $|\sD|^{-k}=(\sD^{2})^{-\frac{n-l}{2}}$ with $l=n-k$. Since $n$ is even the integers $k$ and $l$ have same parity, 
so from~(\ref{eq:NCR.NCR-heat-kernel-asymptotics1}) and~(\ref{eq:NCR.NCR-heat-kernel-asymptotics2}) we get: 
\begin{gather}
    \tr_{\sS}c_{|\sD|^{-k}}(x) =0 \qquad \text{if $k$ is odd}, \\
     2 \tr_{\sS}c_{|\sD|^{-k}}(x)=\Gamma(\frac{k}{2})^{-1}\tr_{\sS} a_{\frac{n-k}{2}}(\sD^{2})(x) \qquad \text{if $k$ is even},
     \label{eq:Even.cD-k.an-k-even}
\end{gather}
where $a_{\frac{n-k}{2}}(\sD^{2})(x)$ is the coefficient of $t^{-\frac{k}{2}}$ in the heat kernel asymptotics~(\ref{eq:NCR.heat-kernel-asymptotics}) 
for~$\sD^{2}$. In particular we see that $ \Vol^{(k)}_{g}M$ vanishes when $n$ is odd. 

On the other hand, the coefficients of its heat kernel asymptotics for $\sD^{2}$ have been thoroughly studied (see, e.g.,~\cite{ABP:OHEIT},  
\cite{BGO:LTHILRSSC}, \cite{Gi:ITHEASIT}). For $j=0,1,\ldots$ we have 
\begin{equation}
  \tr_{\sS}a_{j}(\sD^{2})(x)=(2\pi)^{-\frac{n}{2}}\alpha_{j}(x)dv_{g}(x),
     \label{eq:NCG.heat-form}
\end{equation}
where $\alpha_{j}(x)$ is a local Riemannian invariant of weight $j$ in the sense of~\cite{ABP:OHEIT} and~\cite{Gi:ITHEASIT}. This  means  
that $\alpha_{j}(x)$ is  a linear combination depending only on $n$ of complete tensorial contractions of covariant 
derivatives of the curvature tensor. 

In addition, with the normalization used in~(\ref{eq:NCG.heat-form}) the coefficients of this linear combination don't depend 
on $n$, and in this sense $\alpha_{j}(x)$ does not depend on~$n$. More precisely, by the Lichnerowicz's formula we have
\begin{equation}
   \sD^{2}=(\nabla^{\sS})^{*}\nabla^{\sS}+\frac{1}{4}\kappa, 
%     \label{¥}
\end{equation}
where $\nabla^{\sS}$ denotes the lift to $\cS$ of the Levi-Civita connection. In particular, the operator
$\sD^{2}$ differs from the connection Laplacian $(\nabla^{\sS})^{*}\nabla^{\sS}$ by a scalar section of $\End \sS$. Therefore,  
the functoriality arguments of Gilkey~\cite[Sect.~4.1]{Gi:ITHEASIT} to prove the independence with respect to $n$ of
$(4\pi)^{\frac{n}{2}}a_{j}(\sD^{2})(x)$ apply \emph{verbatim} to prove the independence with respect to $n$ of 
$(4\pi)^{\frac{n}{2}}(\rk \sS)^{-1}\tr_{\sS}a_{j}(\sD^{2})(x)$. Since $\rk \sS=2^{\frac{n}{2}}$ we now see that with the normalization 
of~(\ref{eq:NCG.heat-form}) the local Riemannian invariant $\alpha_{j}(x)$ does not depend on $n$. 

By using the results of~\cite[Sect.~4.1]{Gi:ITHEASIT}, and the references therein, we can get explicit formulas for $\alpha_{j}(x)$ when $j=0,\ldots,5$. 
Let us work out these formulas for $j=0,1,2$.   

First, in the notation of~\cite[Thm.~4.1.6]{Gi:ITHEASIT} we have $E=-\frac{1}{4}\kappa$, so we get  
\begin{equation}
    \alpha_{0}(x)=1 \qquad \text{and} \qquad  \alpha_{1}(x)=\frac{-\kappa(x)}{12}.
    \label{eq:Even.formula-alpha01}
\end{equation}

A formula for $\alpha_{2}(x)$ is given in~\cite{BGO:LTHILRSSC}. Let us explain how to deduce it from the formulas of~\cite[Thm.~4.1.6]{Gi:ITHEASIT}. 
Let $\cR_{ij}=\frac{1}{4}R_{ijkl}c(dx^{k})c(dx^{l})$ be the lift to $\cS$ of the curvature tensor, 
where $c$ is the Clifford action of $\Lambda^{*}_{\C}T^{*}M$ on $\cS$ so that $c(dx^{k})c(dx^{l})+c(dx^{l})c(dx^{k})=-2g^{kl}$. In the notation 
of~\cite[Thm.~4.1.6]{Gi:ITHEASIT} we have $\Omega_{ij}=\cR_{ij}$, so by using~\cite[Thm.~4.1.6]{Gi:ITHEASIT}
we see that $\tr_{\sS}a_{2}(\sD^{2})(x)=(4\pi)^{-\frac{n}{2}}\tr_{\sS}\tilde{\alpha}_{2}(x)dv_{g}(x)$, with
\begin{equation}
     \tilde{\alpha}_{2}(x) = 
    \frac{1}{1440}\left( -12 \Delta_{g}\kappa+5\kappa^{2}-8|\rho|^{2}+8|R|^{2}+120 \cR^{ij}\cR_{ij}\right),
    \label{eq:Even.tilde-alpha2}
\end{equation}
where $|\rho|^{2}=\rho^{jk}\rho_{jk}$  and $|R|^{2}=R^{ijkl}R_{ijkl}$. %are the norm squares of the Ricci and curvature tensors.

We have $\alpha_{2}(x)=2^{-\frac{n}{2}}\tr_{\sS}\tilde{\alpha}_{2}(x)$, so to determine $\alpha_{2}(x)$ we need to evaluate 
$ \tr_{\sS}\cR^{ij}\cR_{ij}$. In normal coordinates centered at $x$ we have
\begin{equation}
     \tr_{\sS}\cR^{ij}\cR_{ij}= \frac{1}{16}\delta_{kr}\delta_{ls} R^{ijrs}R_{ijpq}\tr_{\sS}[c(dx^{k})c(dx^{l})c(dx^{p})c(dx^{q})].     
\end{equation}
For $1\leq i_{1}<\ldots<i_{m}\leq n$ we have $\tr_{\sS}\left[c(dx^{i_{1}})\ldots 
c(dx^{i_{m}})\right]=0$ (see~\cite[Thm.~I.8]{Ge:POSASIT}). Combining this with the Clifford relations 
$c(dx^{k})c(dx^{l})+c(dx^{l})c(dx^{k})=-2\delta^{kl}$ allows us to prove that:
\begin{equation}
    \tr_{\sS}[c(dx^{k})c(dx^{l})c(dx^{p})c(dx^{q})] =\left\{ 
    \begin{array}{cl}
        2^{\frac{n}{2}}& \text{if $(k,p)=(l,q)$ or $(k,l)=(q,p)$},\\
        - 2^{\frac{n}{2}} & \text{if $(k,l)=(p,q)$}, \\
        0 & \text{otherwise}.
    \end{array}\right.
%     \label
\end{equation}
Together with the Bianchi identity $R_{ijlk}=-R_{ijkl}$ this gives 
\begin{equation}
      \tr_{\sS}\cR^{ij}\cR_{ij}= -\frac{2^{\frac{n}{2}}}{8}R^{ijkl}R_{ijkl} = -\frac{2^{\frac{n}{2}}}{8}|R|^{2}.
\end{equation}
Combining this with~(\ref{eq:Even.tilde-alpha2}) and the fact that $ \alpha_{2}(x)=2^{-\frac{n}{2}}\tr_{\sS}\tilde{\alpha}_{2}(x)$ we get
\begin{equation}
    \alpha_{2}(x)= \frac{-1}{1440}\left( 12 \Delta_{g}\kappa^{2}-5\kappa(x)^{2}+8|\rho(x)|^{2}+7|R(x)|^{2}\right).
    \label{eq:Even.formula-alpha2}
\end{equation}

Let us now go back to the computation of $\Vol^{(k)}_{g}M$ when $k$ is even. 
From~(\ref{eq:Even.Volk-cD-k}), (\ref{eq:Even.cD-k.an-k-even}), and~(\ref{eq:NCG.heat-form}) we obtain
\begin{equation}
    \Vol^{(k)}_{g}M= \nu_{n,k}\int_{M}\alpha_{\frac{n-k}{2}}(x)dv_{g}(x), 
    \quad \nu_{n,k}= \frac{(c_{n})^{k}}{n}.\frac{1}{2}\Gamma(\frac{k}{2})^{-1}.(2\pi)^{-\frac{n}{2}}.
     \label{eq:Even.Volkg-nunk}
\end{equation}
Using the expression of $c_{n}$ in~(\ref{eq:NCG.length-element-even}) and the identity $\Gamma(\frac{k}{2}+1)=\frac{k}{2}\Gamma(\frac{k}{2})$ we get
\begin{equation}
    \nu_{n,k}=\frac{1}{n}(2\pi)^{\frac{k}{2}}\Gamma(\frac{n}{2}+1)^{\frac{k}{n}}.k\Gamma(\frac{k}{2}+1)^{-1}.(2\pi)^{-\frac{n}{2}}
    =\frac{k}{n}(2\pi)^{\frac{k-n}{2}}\frac{\Gamma(\frac{n}{2}+1)^{\frac{k}{n}}}{\Gamma(\frac{k}{2}+1)}.
%     \label{¥}
\end{equation}

Summarizing all this we have proved: 

\begin{proposition}\label{prop:Even.main}
1) $\Vol^{(k)}_{g}M$ vanishes when $k$ is odd.\smallskip

2) When $k$ is even we have
\begin{equation}
    \Vol^{(k)}_{g}M=\nu_{n,k}\int_{M}\alpha_{\frac{n-k}{2}}(x)dv_{g}(x), \quad 
    \nu_{n,k}=\frac{k}{n}(2\pi)^{\frac{k-n}{2}}\frac{\Gamma(\frac{n}{2}+1)^{\frac{k}{n}}}{\Gamma(\frac{k}{2}+1)},
%     \label
\end{equation}
where $\alpha_{\frac{n-k}{2}}(x)$ is  a linear combination of complete contractions of
weight $\frac{n-k}{2}$ of covariant derivatives of the curvature tensor. The coefficients of this linear combination depend only on 
$n-k$ and  for $\frac{n-k}{2}=0,1,2$ the invariant $\alpha_{\frac{n-k}{2}}(x)$  
 is explicitly given by~(\ref{eq:Even.formula-alpha01}) and~(\ref{eq:Even.formula-alpha2}). 
\end{proposition}

As a consequence we see that the lower dimensional volumes are integrals of local Riemannian invariants. Since the latter makes sense independently 
of the existence of a spin structure we then can make use of Proposition~\ref{prop:Even.main} to extend the  definition of the lower 
dimensional volumes to \emph{any} compact Riemannian manifold of even dimension as follows. 

\begin{definition}\label{def:Even.general-def}
Let $(M^{n},g)$ be an even-dimensional compact Riemannian manifold.  
Then for $k=1,\ldots,n$ the $k$'th dimensional volume of $(M^{n},g)$ is
\begin{equation}
     \Vol^{(k)}_{g}M:=\left\{
     \begin{array}{ll}
         \nu_{n,k}\int_{M}\alpha_{\frac{n-k}{2}}(x)dv_{g}(x) &  \text{if $k$ is even,}\\
               0 & \text{if $k$ is odd.}\\
     \end{array} \right.
%     \label
\end{equation}
\end{definition}

Notice that this definition is purely differential geometric and does not make appeal anymore to noncommutative geometry. Furthermore, because the 
definition involves the integral of Riemannian invariants we see that two isometric even-dimensional compact Riemannian manifolds have same lower 
dimensional volumes. 

Let us now determine the area  $\op{Area}_{g}M:=\Vol^{(2)}M$ in dimension $4$ and $6$. 
First, we have 
$\nu_{4,2}=\frac{2}{4}(2\pi)^{\frac{2-4}{2}}\frac{\Gamma(\frac{4}{2}+1)^{\frac{2}{4}}}{\Gamma(\frac{2}{2}+1)}= \frac{1}{2\pi \sqrt{2}}$, so by 
using~(\ref{eq:Even.formula-alpha01}) we see that in dimension 4 we have
    \begin{equation}
        \op{Area}_{g}M^{4}=\frac{-1}{24 \pi \sqrt{2}}\int_{M} \kappa(x)dv_{g}(x).
         \label{eq:NCG.area4}
    \end{equation}
This is the formula given by Connes for the area of a 4-dimensional volume that yields a spectral interpretation of the 
Einstein-Hilbert action~(see~\cite{KW:GNGWR}, \cite{Ka:DOG}).

Finally, when $n=6$ we have 
$\nu_{6,2}=\frac{2}{6}(2\pi)^{\frac{2-6}{2}}\frac{\Gamma(\frac{6}{2}+1)^{\frac{2}{6}}}{\Gamma(\frac{2}{2}+1)}= 
    \frac{1}{6}(2\pi)^{-2}\sqrt[3]{6}=\frac{\sqrt[3]{6}}{24\pi^{2}}$.
Together with~(\ref{eq:Even.formula-alpha2}) this shows that in dimension 6 we have 
\begin{equation}
    \op{Area}_{g}M^{6}= \frac{-1}{34560}\int_{M}(-5\kappa(x)^{2}+8|\rho(x)|^{2}+7|R(x)|^{2})dv_{g}(x),
    \label{eq:Even.Area-6}
\end{equation}
where we have used the fact that $\int_{M}\Delta_{g}\kappa dv_{g}(x)=\int_{M}g(\nabla \kappa,\nabla 1)dv_{g}(x)=0$.  

\section{Lower dimensional volumes (odd dimension)}
\label{sec:odd-dimension}
In this section we define the lower dimensional volumes of a Riemannian manifold in odd dimensional. This will follow along similar lines as that of 
the even dimensional case with only few modifications. 

Let $(M^{n},g)$ be a compact Riemannian manifold of odd dimension. As in the even dimensional case we start by assuming that $M$ is spin and we let 
$\sD$ be the Dirac operator acting on the sections of a spin bundle $\sS$ for $M$. Taking into account that in odd dimension we have $\rk \sS=2^{\frac{n-1}{2}}$ 
in the same way as in~(\ref{eq:Even.cD-n-volume-form}) we have
\begin{equation}
    \frac{1}{n}\tr_{\sS} c_{\sD^{-n}}(x)=\frac{(2\pi)^{-n}}{n}\rk \sS.|S^{n-1}|dv_{g}(x)= 
    \frac{(2\pi)^{-\frac{n}{2}}}{\sqrt{2}} \Gamma(\frac{n}{2}+1)^{-1}dv_{g}(x) .
%     \label{¥}
\end{equation}
Therefore, as in~(\ref{eq:NCG.volume-form-|D|-even}) we see that, for any $f \in C^{\infty}(M)$, we have 
\begin{equation}
    \bint f \sD^{-n}=\frac{1}{n}\int_{M} f(x)\tr_{\sS}c_{\sD^{-n}}(x)= 
   \frac{(2\pi)^{-\frac{n}{2}}}{\Gamma(\frac{n}{2}+1)} \int_{M} f(x)dv_{g}(x) .
     \label{eq:Odd.volume-form-|D|}
\end{equation}
Hence the operator $ \sqrt{2}(2\pi)^{\frac{n}{2}}\Gamma(\frac{n}{2}+1) \sD^{-n}$ recaptures the volume form of $M$. 

As in the even dimensional case we define the length element and the lower dimensional volumes of $M$ as follows. 

\begin{definition}
   1) The noncommutative length element of $M$ is 
   \begin{equation}
    ds:=c_{n}'|\sD|^{-1}, \qquad c_{n}'=2^{\frac{1}{2n}}\sqrt{2\pi}\Gamma(\frac{n}{2}+1)^{\frac{1}{n}}.
     \label{eq:Odd.length-element}
\end{equation}
\indent 2)  For $k=1,\ldots,n$ the $k$'th dimensional volume of $M$ is
\begin{equation}
    \Vol^{(k)}_{g}M:=\bint ds^{k}.
%     \label{¥}
\end{equation}
\end{definition}

As in~(\ref{eq:Even.volume-k=n}) using~(\ref{eq:Odd.volume-form-|D|}) we see that $\Vol^{(n)}_{g}M=\Vol_{g}M$. 
In general, we obtain a differential-geometric expression $\Vol^{(k)}_{g}M$ as 
follows.

As in~(\ref{eq:Even.Volk-cD-k}) we have: 
\begin{equation}
    \Vol^{(k)}_{g}M= (c_{n}')^{k}\bint |\sD|^{-k}=  \frac{(c_{n}')^{k}}{n}\int_{M} \tr_{\sS}c_{|\sD|^{-k}}(x).
     \label{eq:Odd.Volk-cD-k}
\end{equation}
Moreover, we have $|\sD|^{-k}=(\sD^{2})^{-\frac{n-l}{2}}$ with $l=n-k$. Here $n$ is odd, so the integers $k$ and $l$ have opposite parities. Therefore, 
from~(\ref{eq:NCR.NCR-heat-kernel-asymptotics1}) and~(\ref{eq:NCR.NCR-heat-kernel-asymptotics2}) we get: 
\begin{gather}
    \tr_{\sS}c_{|\sD|^{-k}}(x) =0 \qquad \text{if $k$ is even}, \\
     2 \tr_{\sS}c_{|\sD|^{-k}}(x)=\Gamma(\frac{k}{2})^{-1}\tr_{\sS} a_{\frac{n-k}{2}}(\sD^{2})(x) \qquad \text{if $k$ is odd}.
     \label{eq:Odd.cD-k.an-k-odd}
\end{gather}
Incidentally, when $k$ is even the $k$'th dimensional volume $ \Vol^{(k)}_{g}M$ is zero. 

Next, in the even dimensional case we explained that the arguments of~\cite[Sect.~4.1]{Gi:ITHEASIT} implies that 
$ (4\pi)^{-\frac{n}{2}}(\rk \sS)^{-1}\tr_{\sS}a_{j}(\sD^{2})(x)$ is a local Riemannian invariant 
independent of $n$, in the sense that in any dimension it is the linear combination with same coefficients of the same complete contractions of 
covariant derivatives of the curvature tensor. 

In fact, Gilkey's arguments go by crossing with $S^{1}$, so $ (4\pi)^{-\frac{n}{2}}(\rk 
\sS)^{-1}\tr_{\sS}a_{j}(\sD^{2})(x)$ is independent of $n$ not only in even dimension, but also in \emph{any} dimension. Therefore  $ (4\pi)^{-\frac{n}{2}}(\rk 
\sS)^{-1}\tr_{\sS}a_{j}(\sD^{2})(x)$ agrees with the local Riemannian invariant $\alpha_{j}(x)$ in~(\ref{eq:NCG.heat-form}). 
Since $\rk \sS=2^{\frac{n-1}{2}}$ it follows that for $j=0,1,,\ldots$ we have
\begin{equation}
    \tr_{\sS}a_{j}(\sD^{2})(x)=\frac{(2\pi)^{-\frac{n}{2}}}{\sqrt{2}}\alpha_{j}(x)dv_{g}(x). 
    \label{eq:Odd.heat-form}
\end{equation}

Now, by combining~(\ref{eq:Odd.Volk-cD-k}), ~(\ref{eq:Odd.cD-k.an-k-odd}) and~(\ref{eq:Odd.heat-form}) we obtain
\begin{equation}
    \Vol^{(k)}_{g}M= \nu'_{n,k}\int_{M}\alpha_{\frac{n-k}{2}}(x)dv_{g}(x), 
    \quad \nu'_{n,k}= \frac{(c_{n}')^{k}}{n}.\frac{1}{2}\Gamma(\frac{k}{2})^{-1}.\frac{(2\pi)^{-\frac{n}{2}}}{\sqrt{2}}.
%     \label{¥}
\end{equation}
As in~(\ref{eq:Even.Volkg-nunk}) using the formula of $c_{n}'$ in~(\ref{eq:Odd.length-element}) we get
\begin{equation}
 \nu'_{n,k}=  \frac{1}{n}2^{\frac{k}{2n}}(2\pi)^{\frac{k}{2}}\Gamma(\frac{n}{2}+1)^{\frac{k}{n}}.k\Gamma(\frac{k}{2}+1)^{-1}. 
    \frac{(2\pi)^{-\frac{n}{2}}}{\sqrt{2}} 
    =\frac{k}{n}2^{a}\pi^{\frac{k-n}{2}}\frac{\Gamma(\frac{n}{2}+1)^{\frac{k}{n}}}{\Gamma(\frac{k}{2}+1)}.
%     \label{¥}
\end{equation}
where $a=\frac{k}{2n}+\frac{k}{2}-\frac{n}{2}-\frac{1}{2}=\frac{(k-n)(n+1)}{2n}$. Therefore, we have proved: 

\begin{proposition}\label{prop:Odd.main}
1)  $\Vol^{(k)}_{g}M$ vanishes when $k$ is even.\smallskip

2) When $k$ is odd we have
\begin{equation}
    \Vol^{(k)}_{g}M=\nu_{n,k}'\int_{M}\alpha_{\frac{n-k}{2}}(x)dv_{g}(x), \quad 
    \nu_{n,k}'=\frac{k}{n}2^{\frac{(k-n)(n+1)}{2n}}\pi^{\frac{k-n}{2}}\frac{\Gamma(\frac{n}{2}+1)^{\frac{k}{n}}}{\Gamma(\frac{k}{2}+1)},
\end{equation}
where $\alpha_{\frac{n-k}{2}}(x)$ is the same local Riemannian invariant as in Proposition~\ref{prop:Even.main}.
\end{proposition}

As with Proposition~\ref{prop:Even.main} this shows that the lower dimensional volumes are integrals of local Riemannian invariants which make sense 
on non-spin manifolds as well. Thus this allows us to extend the definition of the lower dimensional volumes to 
\emph{any} odd-dimensional compact Riemannian manifold as follows.

\begin{definition}\label{def:Odd.general-def}
Let $(M^{n},g)$ be an odd-dimensional compact Riemannian manifold. Then for $k=1,\ldots,n$ the $k$'th dimensional volume of $M$ is: 
\begin{equation}
     \Vol^{(k)}_{g}M:=\left\{
     \begin{array}{ll}
         \nu'_{n,k}\int_{M}\alpha_{\frac{n-k}{2}}(x)dv_{g}(x) &  \text{if $k$ is odd,}\\
               0 & \text{if $k$ is even.}\\
     \end{array} \right.
     \label{eq:Odd.Volkg-non-spin}
\end{equation}
\end{definition}

Let us now give explicit formulas for the length $ \op{Length}_{g}M:=\Vol^{(1)}_{g}M$ in dimension 3 and in dimension 5.

Since $\Gamma(\frac{1}{2}+1)=\frac{1}{2}\Gamma(\frac{1}{2})=\frac{\sqrt{\pi} }{2}$ and 
$\Gamma(\frac{3}{2}+1)=\frac{3}{2}\Gamma(\frac{1}{2}+1)=\frac{3\sqrt{\pi}}{2}$ we have 
\begin{equation}
    \nu'_{3,1}=\frac{1}{3}2^{\frac{(1-3)(3+1)}{2.3}}\pi^{\frac{1-3}{2}}\frac{\Gamma(\frac{3}{2}+1)^{\frac{1}{3}}}{\Gamma(\frac{1}{2}+1)} 
    = \frac{1}{3}2^{-\frac{4}{3}}\pi^{-1}\frac{(3\sqrt{\pi}/2)^{\frac{1}{3}}}{\sqrt{\pi}/2}= \frac{\sqrt[3]{6}}{6\pi\sqrt[3]{\pi}}.
%     \label{¥}
\end{equation}
Combining this with~(\ref{eq:Odd.Volkg-non-spin}) and~(\ref{eq:Even.formula-alpha01}) then shows that in dimension 3 we have
\begin{equation}
    \op{Length}_{g}M^{3}=-\frac{\sqrt[3]{6}}{72\pi\sqrt[3]{\pi}}\int_{M}\kappa(x)dv_{g}(x).
    %     \label{¥}
\end{equation}

Finally, since $\Gamma(\frac{5}{2}+1)=\frac{5}{2}\Gamma(\frac{3}{2}+1)=\frac{15\sqrt{\pi}}{4}$ we get
\begin{equation}
     \nu'_{5,1}=\frac{1}{5}2^{\frac{(1-5)(5+1)}{2.5}}\pi^{\frac{1-5}{2}}\frac{\Gamma(\frac{5}{2}+1)^{\frac{1}{5}}}{\Gamma(\frac{1}{2}+1)} 
     = \frac{1}{5}2^{-\frac{12}{5}}\pi^{-2}\frac{(15\sqrt{\pi}/4)^{\frac{1}{5}}}{\sqrt{\pi}/2}= \frac{\pi\sqrt[5]{30}}{20\sqrt[10]{\pi}}.
 %     \label{¥}
\end{equation}
Therefore, by using~(\ref{eq:Odd.Volkg-non-spin}) and~(\ref{eq:Even.formula-alpha2}) we see that in dimension 5 we have
\begin{equation}
     \op{Length}_{g}M^{5}= -\frac{\pi\sqrt[5]{60}}{28800\sqrt[10]{\pi}} \int_{M}(-5\kappa(x)^{2}+8|\rho(x)|^{2}+7|R(x)|^{2})dv_{g}(x),
%     \label{¥}
\end{equation}
where as in~(\ref{eq:Even.Area-6}) we have used the fact that $\int_{M}\Delta_{g}\kappa(x) dv_{g}(x)=0$.

\end{document}